\documentclass[reqno]{amsart}
\usepackage{hyperref}
\usepackage[russian]{babel}
\usepackage{amssymb}

\begin{document}

\phantom{A}

\title[On optimal control in a model of rigid-viscoplastic media]
{On optimal control in a model \\ of rigid-viscoplastic media \\with Dirichlet boundary conditions}

\author[M. A. Artemov, A. V. Skobaneva]
{Mikhail A. Artemov, Anna V. Skobaneva}

\address{Mikhail A. Artemov \newline
Department of Applied Mathematics, Informatics and Mechanics,Voronezh State University, Voronezh,  Russia}
\email{artemov\_m\_a@mail.ru}

\address{Anna V. Skobaneva \newline
Department of Applied Mathematics, Informatics and Mechanics,Voronezh State University, Voronezh,  Russia}

\thanks{Submitted August 22, 2017.}

\maketitle
\newtheorem{lemma}{Лемма}
\newtheorem{theorem}{Теорема}
\newtheorem{definition}{Определение}
\newtheorem{corollary}{Следствие}
\newtheorem{proposition}{Предложение}
\allowdisplaybreaks

\begin{quote}
{\bf Abstract.} In this paper, we consider the optimal control problem in a 3D flow model for incompressible rigid-viscoplastic media of the Bingham kind with homogeneous Dirichlet boundary conditions and a given cost functional. On the basis of methods of the theory of variational inequalities with pseudo\-monotone operators, a theorem on the solvability of the optimization problem in the class of weak steady solutions is proved.
\end{quote}
\vspace{4mm}

\begin{quote}
{\bf Аннотация.}  В статье рассмотрена задача оптимального управления в трехмерной модели течения несжимаемой жестко-вязко-пластической среды типа Бингама с однородными краевыми условиями Дирихле и заданным целевым функционалом. На основе методов теории вариационных неравенств с псевдомонотонными операторами доказана теорема о разрешимости задачи оптимизации в классе <<слабых>> стационарных решений.
\end{quote}

\section{Введение}
Жестко-вязко-пластические среды типа Бингама \cite{PRA, DUV} характеризуются тем, что в них при  малых напряжениях тензор скоростей деформаций равен {\bf 0} и в~соответствующих областях среда движется как твердое тело до тех пор, пока некоторая функция напряжений не достигнет своего предела текучести; выше этого предела среда ведет себя уже  как несжимаемая вязкая жидкость.

 Модели бингамовских жидкостей достаточно универсальны и уже долгое время успешно применяются при описании  потоков большого числа реальных вязкопластических сред, в том числе таких материалов, как цементы, бетон, суспензии, пасты, краски, гели, некоторые виды нефтей и масел.

Математическое изучение уравнений динамики сред Бингама началось в работах французских ученых Г. Дюво и Ж.-Л. Лионса. Полученные ими результаты подробно изложены в классической монографии  \cite{DUV}. В настоящий момент исследования по данному направлению активно продолжаются; интерес к таким задачам стимулируется разнообразными приложениями как в математической гидродинамике, так в и ряде технологических процессов. Из последних опубликованных работ, связанных с вышеупомянутой моделью и ее обобщениями, можно отметить \cite{MAL, MAM, LAP, KHO, SHE}.

В предлагаемой статье ставится задача об оптимальном управлении  трехмерным  стационарным течением жестко-вязко-пластической среды  типа Бингама в ограниченной области ${{\mathbb{R}}}^3$ при условии прилипания на границе области. Используя методы теории вариационных неравенств с псевдомонотонными операторами, мы доказываем теорему о разрешимости задачи оптимизации в классе <<слабых>> решений.

Данная работа развивает результаты, полученные в \cite{BAR4} для случая нелинейно-вязких жидкостей. Отметим также статьи \cite{WAC, BAR1, BAR2, BAR3}, в которых рассматриваются близкие по постановке задачи оптимизации для моделей неньютоновских сред.

\section{Формулировка задачи оптимального управления} 
Предположим, что $\Omega$ --- ограниченное открытое множество  в ${{\mathbb{R}}}^3$ с достаточно гладкой границей $\partial \Omega$. Будем рассматривать
задачу оптимизации стационарных трехмерных потоков несжимаемой жестко-вязко-пластической среды типа Бингама в области $\Omega$:
\begin{equation}
\sum_{i=1}^3 v_i\frac{\partial v_j}{\partial x_i}-(\text{div}{\boldsymbol{\sigma}})_j+\frac{\partial p}{\partial x_j}={f_j}+{u_j},\quad j\in\{1,2,3\}, \label{bing-op-1}
\end{equation} 
\begin{equation}
\text{div}\,\boldsymbol{v}=\sum_{i=1}^3\frac{\partial v_i}{\partial x_i}=0, \label{bing-op-2}
\end{equation} 
\begin{equation}
{\sigma}_{ij}=\mu(|\boldsymbol{\mathcal{E}}|){\mathcal{E}}_{ij}+g\frac{{\mathcal{E}}_{ij}}{|\boldsymbol{\mathcal{E}}|},\mbox{ если }|\boldsymbol{\mathcal{E}}|\neq 0,\quad i,j\in\{1,2,3\}, \label{bing-op-3}
\end{equation} 
\begin{equation}
|\boldsymbol{\sigma}|\leq g,\mbox{ если }|\boldsymbol{\mathcal{E}}|=0, \label{bing-op-4}
\end{equation} 
\begin{equation}
\boldsymbol{v}|_{\partial \Omega}=\boldsymbol{0},\label{bing-op-5}
\end{equation} 
\begin{equation}
\boldsymbol{u}\in\boldsymbol{U}_{\rm ad},\label{bing-op-6}
\end{equation} 
\begin{equation}
J(\boldsymbol{v},\boldsymbol{u})\to \min.\label{bing-op-7}
\end{equation} 

В системе \eqref{bing-op-1}--\eqref{bing-op-7} используются обозначения: 
$\boldsymbol{v}=(v_1(\boldsymbol{x}),v_2(\boldsymbol{x}),v_3(\boldsymbol{x}))$ --- скорость течения в точке $\boldsymbol{x}=(x_1, x_2, x_3)$ из области $\Omega$, $p=p(\boldsymbol{x})$ --- давление,  $\boldsymbol{\sigma}=({\sigma}_{ij}(\boldsymbol{x}))$ --- девиатор тензора напряжений, $\boldsymbol{f}=(f_1(\boldsymbol{x}),f_2(\boldsymbol{x}),f_3(\boldsymbol{x}))$ --- внешняя сила,  действующая на среду,   $\text{div}\,\boldsymbol{\sigma}$ --- вектор с компонентами
$$
(\text{div}\,\boldsymbol{\sigma})_j=\sum_{i=1}^3\limits\frac{\partial \sigma_{ij}}{\partial x_i},
$$
$\boldsymbol{\mathcal{E}}=\boldsymbol{\mathcal{E}}(\boldsymbol{v})$ --- тензор скоростей деформации, 
$$
{\mathcal{E}}_{ij}={\mathcal{E}}_{ij}(\boldsymbol{v})=\frac{1}{2}\left(\frac{\partial v_i}{\partial x_j}+\frac{\partial v_j}{\partial x_i}\right),
$$
$|\boldsymbol{\mathcal{E}}|^2$ --- второй инвариант тензора $\boldsymbol{\mathcal{E}}$, определяемый по формуле
\begin{equation*}
|\boldsymbol{\mathcal{E}}|^2=\boldsymbol{\mathcal{E}}:\boldsymbol{\mathcal{E}}=\sum_{i,j=1}^3\mathcal{E}_{ij}^2,
\end{equation*} 
$\mu(|\boldsymbol{\mathcal{E}}|)>0$ --- вязкость, $g=g(\boldsymbol{x})>0$ --- порог текучести, разделяющий два типа поведения среды,
$\boldsymbol{u}=(u_1(\boldsymbol{x}),u_2(\boldsymbol{x}),u_3(\boldsymbol{x}))$ --- набор управляющих параметров, 
$\boldsymbol{U}_{\rm ad}$ --- множество допустимых управлений, $J$ --- целевой функционал.

При рассмотрении задачи \eqref{bing-op-1}--\eqref{bing-op-7} мы считаем, что $\boldsymbol{f}$, $\mu$, $g$, $\boldsymbol{U}_{\rm ad}$, $J$ --- известны, а  распределение скоростей $\boldsymbol{v}$ и управление $\boldsymbol{u}$ (вместе с ассоциированным с ними давлением~$p$) являются неизвестными величинами.
\vspace{0.7mm}

{\bf Замечание 1.} Данная задача может быть отнесена к так называемым {\it задачам со свободной границей}, т.е. таким задачам, в которых уравнения различны в различных частях рассматриваемой области. При этом неизвестная заранее свободная граница отделяет зону области $\Omega$, где материал жесткий, от зоны, где материал проявляет жидкие свойства. Свободная граница участвует в постановке задачи неявно и может быть определена только после нахождения решений.
\vspace{0.7mm}

{\bf Замечание 2.} Если положить $g\equiv 0$, то рассматриваемая реологическая модель сводится к модели нелинейно-вязкой жидкости~\cite{LIT}; если к тому же  $\mu={\rm const}$, то мы имеем дело с хорошо известной системой уравнений Навье--Стокса \cite{TEM}, описывающей динамику классической ньютоновской жидкости. При возраствании предела текучести $g$ в потоке появляются области, где жидкость ведет себя подобно твердому телу, а при значительном увеличении $g$ течение полностью блокируется.

\section{Используемые функциональные пространства} 
Далее по тексту будут использоваться пространства Лебега ${L}_r(\Omega)$ и Соболева ${W}^{m}_{q}(\Omega)$, состоящие из функций, заданных в области $\Omega$.  Нормы в этих пространствах определяются обычным образом (см.\cite{TEM, ADA}). Когда речь идет о пространствах векторнозначных функций, будем использовать жирный шрифт: $\boldsymbol{L}_r(\Omega)$, $\boldsymbol{W}^{m}_{q}(\Omega)$ и т.д.

Пусть $\mathcal{D}(\Omega)$ --- множество функций класса $\mathcal{C}^\infty$ с компактным носителем, содержащимся в области $\Omega$. Замыкание множества $\mathcal{D}(\Omega)$ в ${W}^{1}_{ 2}(\Omega)$ обозначим через ${H}^1_0(\Omega)$. 

Из неравенства Пуанкаре \cite{TEM} 
$$
\|{v}\|_{{L}_2(\Omega)}\leq C_1 \|\nabla{v}\|_{\boldsymbol{L}_2(\Omega)}\quad C_1={\rm const},\;\forall\,{v}\in{H}^1_0(\Omega)
$$
следует, что в ${H}^1_0(\Omega)$ можно ввести норму
$$
\|{v}\|_{{H}^1_0(\Omega)}= \|\nabla{v}\|_{\boldsymbol{L}_2(\Omega)}
$$
и эта норма будет эквивалента норме $\|\bullet\|_{{W}^1_2(\Omega)}$.

Введем также в рассмотрение также пространство соленоидальных функций
$$
\boldsymbol{V}=\{\boldsymbol{v}\in \boldsymbol{H}^1_0(\Omega):\;\;\text{div}\,\boldsymbol{v}=0\}
$$
со скалярным произведением 
$$
(\boldsymbol{v}, \boldsymbol{w})_{\boldsymbol{V}}=\int\limits_\Omega\boldsymbol{\mathcal{E}}(\boldsymbol{v}):\boldsymbol{\mathcal{E}}(\boldsymbol{w})\,dx
$$
и нормой
$
\|\boldsymbol{v}\|_{\boldsymbol{V}}=(\boldsymbol{v}, \boldsymbol{v})_{\boldsymbol{V}}^{1/2}.
$
 Из неравенства Корна \cite{DUV} для пространства $\boldsymbol{H}^1_0(\Omega)$
$$
\|\boldsymbol{v}\|_{\boldsymbol{H}^1_0(\Omega)}\leq C_2\|\boldsymbol{\mathcal{E}}(\boldsymbol{v})\|_{\boldsymbol{L}_2(\Omega)}\quad C_2={\rm const},\;\forall\,\boldsymbol{v}\in\boldsymbol{H}^1_0(\Omega),
$$
имеющего фундаментальное значение во многих задачах теории пластичности и вязко-упругости,  вытекает, что введенная выше норма  $\|\bullet\|_{\boldsymbol{V}}$  эквивалентна норме $\|\bullet\|_{\boldsymbol{H}^1_0(\Omega)}$. 

Через $\mathbb{M}^{3\times 3}_{\rm sym}$ обозначим пространство симметрических $3\times 3$-матриц со скалярным произведением
$$
(\boldsymbol{X},\boldsymbol{Y})_{\mathbb{M}^{3\times 3}_{\rm sym}}=\boldsymbol{X}:\boldsymbol{Y}=\sum_{i,j=1}^3X_{ij}Y_{ij}
$$
 и евклидовой нормой 
$$
|\boldsymbol{X}|=(\boldsymbol{X},\boldsymbol{X})_{\mathbb{M}^{3\times 3}_{\rm sym}}^{1/2}.
$$

Как обычно, символы $\to$ и $\rightharpoonup$ обозначают соответственно сильную и слабую сходимость.

 \section{Основные предположения}
Опишем теперь основные наши предположения относительно <<данных>> задачи \eqref{bing-op-1}--\eqref{bing-op-7}. Далее считаем, что
\begin{itemize}
\item[(i)] {\it множество допустимых управлений $\boldsymbol{U}_{\rm ad}\neq\varnothing$ ограничено и секвенциально слабо замкнуто в $\boldsymbol{L}_2(\Omega)$}; 
\item[(ii)] {\it  целевой функционал $J\colon\boldsymbol{V}\times\boldsymbol{L}_2(\Omega)\to{{\mathbb{R}}}$ полунепрерывен снизу относительно слабой сходимости в $\boldsymbol{V}\times\boldsymbol{L}_2(\Omega)$, иными словами:
для любой последовательности $(\boldsymbol{v}_n,\boldsymbol{u}_n)$ такой, что $(\boldsymbol{v}_n,\boldsymbol{u}_n)\rightharpoonup  (\boldsymbol{v},\boldsymbol{u})$ в $\boldsymbol{V}\times\boldsymbol{L}_2(\Omega)$ при $n\to\infty$, имеет место оценка}
\begin{equation*}
J(\boldsymbol{v},\boldsymbol{u})\leq \liminf_{n\to\infty}J(\boldsymbol{v}_n,\boldsymbol{u}_n);
\end{equation*}
\item[(iii)] {\it  выполнено неравенство}
$$
\bigl(\mu(|\boldsymbol{X}|)\boldsymbol{X}-\mu(|\boldsymbol{Y}|)\boldsymbol{Y}\bigr):\bigl(\boldsymbol{X}-\boldsymbol{Y}\bigr)\geq 0\qquad\forall\,\boldsymbol{X},\boldsymbol{Y}\in\mathbb{M}^{3\times 3}_{\rm sym};
$$
\item[(iv)] {\it функция $\mu$ измерима, и существуют  константы $\mu_0$ и $\mu_1$ такие, что }
$$
0<\mu_0<\mu(s)<\mu_1,\quad\forall\,s\in{{\mathbb{R}}}_+;
$$
\item[(v)] {\it выполнено неравенство $g(\boldsymbol{x})>0$ для п.в. $\boldsymbol{x}\in \Omega$ и $g\in L_2(\Omega)$};
\item[(vi)] {\it имеет место включение $\boldsymbol{f}\in\boldsymbol{L}_2(\Omega)$.}
\end{itemize}
\vspace{1mm}

{\bf Замечание 3.} Для выполнения условия (iii) достаточно потребовать, чтобы функция $\mu$ была неубывающей на ${{\mathbb{R}}}_+$, например, $\mu(s)\equiv\arctan(s)+\mu_0$. В самом деле, с применением неравенства Коши--Буняковского--Шварца нетрудно получить, что
$$
\bigl(\mu(|\boldsymbol{X}|)\boldsymbol{X}-\mu(|\boldsymbol{Y}|)\boldsymbol{Y}\bigr):\bigl(\boldsymbol{X}-\boldsymbol{Y}\bigr)
$$
$$
=\mu(|\boldsymbol{X}|)|\boldsymbol{X}|^2-\mu(|\boldsymbol{X}|)\boldsymbol{X}:\boldsymbol{Y}-\mu(|\boldsymbol{Y}|)\boldsymbol{X}:\boldsymbol{Y}+\mu(|\boldsymbol{Y}|)|\boldsymbol{Y}|^2
$$
$$
\geq\mu(|\boldsymbol{X}|)|\boldsymbol{X}|^2-\mu(|\boldsymbol{X}|)|\boldsymbol{X}||\boldsymbol{Y}|-\mu(|\boldsymbol{Y}|)|\boldsymbol{X}||\boldsymbol{Y}|+\mu(|\boldsymbol{Y}|)|\boldsymbol{Y}|^2
$$
$$
=\bigl\{\mu(|\boldsymbol{X}|)|\boldsymbol{X}|-\mu(|\boldsymbol{Y}|)|\boldsymbol{Y}|\bigr\}(|\boldsymbol{X}|-|\boldsymbol{Y}|)\geq  0 \quad\forall\,\boldsymbol{X},\boldsymbol{Y}\in\mathbb{M}^{3\times 3}_{\rm sym},
$$
если функция $\mu$ является неубывающей на ${{\mathbb{R}}}_+$.
\vspace{1mm}

{\bf Замечание 4.}  Типичным примером целевого функционала, удовлетворяющего условию (ii), служит функционал
$$
J(\boldsymbol{v},\boldsymbol{u})=\lambda_1\|\boldsymbol{v}-\widetilde{\boldsymbol{v}}\|_{\boldsymbol{V}}^2+\lambda_2\|\boldsymbol{u}-\widetilde{\boldsymbol{u}}\|_{\boldsymbol{L}_2(\Omega)}^2,
$$
где $\widetilde{\boldsymbol{v}}$ и  $\widetilde{\boldsymbol{u}}$ --- заданные векторнозначные функции, $\lambda_1$ и $\lambda_2$ --- неотрицательные параметры.

\section{Вариационная формулировка задачи:\\ допустимые пары и~оптимальные решения}
 Введем теперь понятие допустимой пары <<скорость-управление>>. 

Пусть
$$
\phi_g\colon\boldsymbol{V}\to {{\mathbb{R}}_+},\quad\phi_g(\boldsymbol{v})=\int\limits_\Omega g|\boldsymbol{\mathcal{E}}(\boldsymbol{v})|\, dx.
$$
\begin{definition}
{\rm Под} допустимой парой {\rm  задачи оптимизации~\eqref{bing-op-1}--\eqref{bing-op-7} понимается
 пара векторнозначных функций $(\boldsymbol{v}, \boldsymbol{u})\in \boldsymbol{V}\times \boldsymbol{U}_{\rm ad}$ такая, что
\begin{gather}
-\sum_{i=1}^3\int\limits_\Omega{v}_i\boldsymbol{v}\cdot \frac{\partial\boldsymbol{w}}{\partial x_i}\,dx+\int\limits_\Omega\mu(|\boldsymbol{\mathcal{E}}(\boldsymbol{v})|)\boldsymbol{\mathcal{E}}(\boldsymbol{v}):\boldsymbol{\mathcal{E}}(\boldsymbol{w}-\boldsymbol{v})\,dx\nonumber
\\
+\phi_g(\boldsymbol{w})-\phi_g(\boldsymbol{v})
\geq \int\limits_\Omega(\boldsymbol{f}+\boldsymbol{u})\cdot(\boldsymbol{w}-\boldsymbol{v})\,dx\label{wks}
\end{gather}
для любой пробной векторнозначной функции $\boldsymbol{w}\in \boldsymbol{V}$.}
\end{definition}

Вариационное неравенство \eqref{wks} определяет так называемые {\it слабые обобщенные решения} задачи.  По поводу корректности данного определения отсылаем читателя к монографии [2].

Пусть  $\boldsymbol{G}$  --- множество всех допустимых пар задачи \eqref{bing-op-1}--\eqref{bing-op-7}.

\begin{definition}
Оптимальным решением {\rm задачи  \eqref{bing-op-1}--\eqref{bing-op-7} назовем пару вектор-функций $(\boldsymbol{v}_0, \boldsymbol{u}_0)\in\boldsymbol{G}$, которая характеризуется равенством}
$$
J(\boldsymbol{v}_0, \boldsymbol{u}_0)=\inf_{(\boldsymbol{v},\boldsymbol{u})\in \boldsymbol{G}}J(\boldsymbol{v},\boldsymbol{u}).
$$
\end{definition}

\section{Теорема о разрешимости задачи оптимального управления} 

Основной результат работы формулируется в виде следующей теоремы.

\begin{theorem}\label{teor1}
При выполнении приведенных выше условий {\rm (i)--(vi)} задача  \eqref{bing-op-1}--\eqref{bing-op-7} имеет хотя бы одно оптимальное решение $(\widehat{\boldsymbol{v}},\widehat{\boldsymbol{u}})\in\boldsymbol{G}$ и 
\begin{equation}
\int\limits_\Omega \mu(|\boldsymbol{\mathcal{E}}(\widehat{\boldsymbol{v}})|)\:|\boldsymbol{\mathcal{E}}(\widehat{\boldsymbol{v}})|^2\,dx+\int\limits_\Omega g|\boldsymbol{\mathcal{E}}(\widehat{\boldsymbol{v}})|\,dx= \int\limits_\Omega(\boldsymbol{f}+\widehat{\boldsymbol{u}})\cdot\widehat{\boldsymbol{v}}\,dx.\label{est-gr}
\end{equation}
\end{theorem}

Доказательство данного результата базируется на теореме о разрешимости вариационных неравенств с псевдомонотонными операторами и обобщенной теореме Вейерштрасса. Для удобства читателя приведем формулировки этих утверждений.

\begin{lemma}[см. \cite{LIO}] \label{lemm1}
Пусть $\boldsymbol{E}$ --- рефлексивное B-пространство, $\boldsymbol{E}^\prime$ --- пространство, сопряженное к $\boldsymbol{E}$, ${\bf A}\colon\boldsymbol{E}\to\boldsymbol{E}^\prime$ --- псевдомонотонный оператор, а $\varphi\colon\boldsymbol{E}\to{{\mathbb{R}}}$ --- выпуклая полунепрерывная снизу функция. Будем предполагать, что имеет место сходимость
$$
\frac{\langle {\bf A}(\boldsymbol{w}),\boldsymbol{w}\rangle+\varphi(\boldsymbol{w})}
{\|\boldsymbol{w}\|_{\boldsymbol{E}}}\to+\infty,
$$
когда $\|\boldsymbol{w}\|_{\boldsymbol{E}}\to+\infty$. 
Тогда для заданного элемента $\boldsymbol{h}\in \boldsymbol{E}^\prime$ существует решение $\boldsymbol{v}\in\boldsymbol{E}$ вариационного неравенства
$$
\langle{\bf A}(\boldsymbol{v})-\boldsymbol{h},\boldsymbol{w}-\boldsymbol{v}\rangle+\varphi(\boldsymbol{w})-\varphi(\boldsymbol{v})\geq 0\qquad\forall\boldsymbol{w}\in\boldsymbol{E}.
$$
\end{lemma}

\begin{lemma}[обобщенная теорема Вейерштрасса \cite{ZEI}] \label{lemm2}
Пусть $\boldsymbol{E}$ --- рефлексивное B-пространство, $\boldsymbol{Q}\subset\boldsymbol{E}$ --- ограниченное и секвенциально слабо замкнутое множество. Предположим, что функционал $\mathcal{J}:\boldsymbol{Q}\to{{\mathbb{R}}}$ полунепрерывен снизу относительно слабой сходимости в $\boldsymbol{E}$. 
Тогда существует элемент $\boldsymbol{y}_0\in\boldsymbol{Q}$ такой, что
$$
\mathcal{J}(\boldsymbol{y}_0)=\inf_{\boldsymbol{y}\in \boldsymbol{Q}}\mathcal{J}(\boldsymbol{y}).
$$
\end{lemma}

 {\bf Доказательство теоремы~\ref{teor1}.}

{\bf Шаг 1.} Докажем сначала, что множество допустимых пар $\boldsymbol{G}$ непусто. Зафиксируем (временно) некоторый элемент $\boldsymbol{u}\in\boldsymbol{U}_{\rm ad}$ и введем в рассмотрение два нелинейных оператора:
$$
\mathbf{M}\colon\boldsymbol{V}\to \boldsymbol{V}^\prime,\quad\langle\mathbf{M}(\boldsymbol{v}),\boldsymbol{w}\rangle=\int\limits_\Omega\mu(|\boldsymbol{\mathcal{E}}(\boldsymbol{v})|)\boldsymbol{\mathcal{E}}(\boldsymbol{v}):\boldsymbol{\mathcal{E}}(\boldsymbol{w})\,dx,
$$
$$
\mathbf{T}_{\boldsymbol{u}}\colon\boldsymbol{V}\to \boldsymbol{V}^\prime,\quad\langle\mathbf{T}_{\boldsymbol{u}}(\boldsymbol{v}),\boldsymbol{w}\rangle=-\sum_{i=1}^3\int\limits_\Omega{v}_i\boldsymbol{v}\cdot \frac{\partial\boldsymbol{w}}{\partial x_i}\,dx-\int\limits_\Omega(\boldsymbol{f}+\boldsymbol{u})\cdot\boldsymbol{w}\,dx.
$$

С помощью этих операторов неравенство \eqref{wks} можно переписать в виде:
\begin{equation}
\langle\mathbf{M}(\boldsymbol{v})+\mathbf{T}_{\boldsymbol{u}}(\boldsymbol{v}),\boldsymbol{w}-\boldsymbol{v}\rangle+\phi_g(\boldsymbol{w})-\phi_g(\boldsymbol{v})\geq 0\quad \forall\,\boldsymbol{w}\in\boldsymbol{V}.\label{eqq-1}
\end{equation}

В силу условия (iii) оператор $\mathbf{M}$  обладает свойством монотонности. Кроме того, данный оператор семинепрерывен. Отсюда вытекает (см. \cite{LIO}), что $\mathbf{M}$ является псевдомонотонным оператором. 

Заметим также, что  $\mathbf{T}_{\boldsymbol{u}}$ --- усиленно непрерывный оператор, в чем несложно убедиться, принимая во внимание компактность вложения Соболева $\boldsymbol{W}^{1}_{2}(\Omega)\hookrightarrow\boldsymbol{L}_4(\Omega)$ (см. \cite{ADA}). Поэтому $\mathbf{M}+\mathbf{T}_{\boldsymbol{u}}$ принадлежит классу псевдомонотонных операторов как сумма псевдомонотонного и усиленно непрерывного операторов (см. \cite{LIO}).

Используя условие (iv) и равенство
$$
\sum_{i=1}^3\int\limits_\Omega{v}_i\boldsymbol{v}\cdot \frac{\partial\boldsymbol{v}}{\partial x_i}\,dx=0,
$$
получаем, что
$$
\frac{\langle {\bf M}(\boldsymbol{v})+\mathbf{T}_{\boldsymbol{u}}(\boldsymbol{v}),\boldsymbol{v}\rangle+\phi_g(\boldsymbol{v})}
{\|\boldsymbol{v}\|_{\boldsymbol{V}}}\to+\infty,
$$
когда $\|\boldsymbol{v}\|_{\boldsymbol{V}}\to+\infty$. 

Применяя лемму~\ref{lemm1}, мы делаем вывод о том, что неравенство \eqref{eqq-1} имеет в пространстве $\boldsymbol{V}$ одно решение или несколько решений, одно из которых обозначим через $\boldsymbol{v}_{\boldsymbol{u}}$. Очевидно, что $(\boldsymbol{v}_{\boldsymbol{u}}, \boldsymbol{u})\in \boldsymbol{G}$.

{\bf Шаг 2.} Покажем, что множество $\boldsymbol{G}$ ограничено в пространстве $\boldsymbol{V}\times\boldsymbol{L}_2(\Omega)$. Возьмем произвольный элемент $(\boldsymbol{v}, \boldsymbol{u})\in \boldsymbol{G}$. В силу наших предположений множество $\boldsymbol{U}_{\rm ad}$ ограничено в $\boldsymbol{L}_2(\Omega)$ (см. условие (i)), поэтому нам необходимо оценить только норму скорости $\boldsymbol{v}$ в пространстве $\boldsymbol{V}$. Для этого полагаем в неравенстве \eqref{wks} $\boldsymbol{w}=\boldsymbol{0}$, а затем $\boldsymbol{w}=-\boldsymbol{v}$; из полученных соотношений выводим
\begin{equation}
\int\limits_\Omega \mu(|\boldsymbol{\mathcal{E}}({\boldsymbol{v}})|)\:|\boldsymbol{\mathcal{E}}({\boldsymbol{v}})|^2\,dx+\int\limits_\Omega g|\boldsymbol{\mathcal{E}}({\boldsymbol{v}})|\,dx= \int\limits_\Omega(\boldsymbol{f}+{\boldsymbol{u}})\cdot{\boldsymbol{v}}\,dx,\label{eqq-20}
\end{equation}
откуда, принимая в расчет условие (iv), выводим неравенство
$$
 \mu_0\|\boldsymbol{v}\|_{\boldsymbol{V}}^2\leq C_\Omega \|\boldsymbol{f}+{\boldsymbol{u}}\|_{\boldsymbol{L}_2(\Omega)}\|\boldsymbol{v}\|_{\boldsymbol{V}},\quad C_\Omega={\rm const}.
$$
Таким образом, мы получаем оценку нормы $\boldsymbol{v}$ в терминах ${\Omega,\boldsymbol{U}_{\rm ad},\boldsymbol{f},\mu_0}$:
$$
\|\boldsymbol{v}\|_{\boldsymbol{V}}\leq \frac{C_\Omega}{\mu_0}\sup_{\boldsymbol{u}\in\boldsymbol{U}_{\rm ad}} \|\boldsymbol{f}+\boldsymbol{u}\|_{\boldsymbol{L}_2(\Omega)}\leq C({\Omega,\boldsymbol{U}_{\rm ad},\boldsymbol{f},\mu_0}).
$$
Тем самым доказана ограниченность $\boldsymbol{G}$.

{\bf Шаг 3.} Покажем теперь, что множество $\boldsymbol{G}$ секвенциально слабо замкнуто в пространстве $\boldsymbol{V}\times\boldsymbol{L}_2(\Omega)$. Возьмем произвольную последовательность $(\boldsymbol{v}_n,\boldsymbol{u}_n)\in\boldsymbol{G}$ такую, что $(\boldsymbol{v}_n,\boldsymbol{u}_n)\rightharpoonup (\boldsymbol{v}_*,\boldsymbol{u}_*)$ в $\boldsymbol{V}\times\boldsymbol{L}_2(\Omega)$ при $n\to \infty$, и проверим, что $(\boldsymbol{v}_*,\boldsymbol{u}_*)\in\boldsymbol{G}$.

Мы, очевидно, имеем слабую сходимость: $\boldsymbol{v}_n\rightharpoonup \boldsymbol{v}_*$ в $\boldsymbol{V}$ и $\boldsymbol{u}_n\rightharpoonup \boldsymbol{u}_*$ в $\boldsymbol{L}_2(\Omega)$ и, более того, благодаря теореме о компактности вложения Соболева $\boldsymbol{W}^{1}_{2}(\Omega)\hookrightarrow\boldsymbol{L}_4(\Omega)$ имеет место также сильная сходимость: $\boldsymbol{v}_n\to \boldsymbol{v}_*$ в $\boldsymbol{L}_4(\Omega)$ при $n\to\infty$.

Ввиду (i) получаем, что $\boldsymbol{u}_*\in\boldsymbol{U}_{\rm ad}$. Далее, поскольку $(\boldsymbol{v}_n,\boldsymbol{u}_n)\in\boldsymbol{G}$ для любого $n\in\mathbb{N}$, то выполнено неравенство:
\begin{equation}
\langle\mathbf{M}(\boldsymbol{v}_n)+\mathbf{T}_{\boldsymbol{u}_n}(\boldsymbol{v}_n),\boldsymbol{w}-\boldsymbol{v}_n\rangle+\phi_g(\boldsymbol{w})-\phi_g(\boldsymbol{v}_n)\geq 0\quad \forall\,\boldsymbol{w}\in\boldsymbol{V},\;n\in\mathbb{N}.\label{eqq-201}
\end{equation}

Нам нужно выполнить предельный переход в \eqref{eqq-201}, устремляя $n\to\infty$.
Для этого заметим следующее. 

Введем B-пространство $\boldsymbol{V}_g$, состоящее из элементов пространства $\boldsymbol{V}$, с нормой
$$
\|\boldsymbol{v}\|_{\boldsymbol{V}_g}=\int\limits_\Omega g|\boldsymbol{\mathcal{E}}(\boldsymbol{v})|\, dx.
$$
Аксиомы нормы выполнены благодаря условию (v). К тому же справедливо соотношение 
$$
\|\boldsymbol{v}\|_{\boldsymbol{V}_g}\leq\|g\|_{L_2(\Omega)}\|\boldsymbol{v}\|_{\boldsymbol{V}}\quad \forall\boldsymbol{v}\in\boldsymbol{V},
$$
т. е. имеет место непрерывное вложение $\boldsymbol{V}\hookrightarrow\boldsymbol{V}_g$. 

Отсюда делаем вывод о том, что $\boldsymbol{v}_n\rightharpoonup \boldsymbol{v}_*$ в $\boldsymbol{V}_g$ при $n\to \infty$ и 
$$
\|\boldsymbol{v}_*\|_{\boldsymbol{V}_g}\leq\liminf_{n\to \infty}\|\boldsymbol{v}_n\|_{\boldsymbol{V}_g},
$$
или, другими словами, 
\begin{equation}
\phi_g(\boldsymbol{v}_*)\leq\liminf_{n\to \infty}\phi_g(\boldsymbol{v}_n).\label{eqq-401}
\end{equation}

Далее, полагаем  $\boldsymbol{w}=\boldsymbol{v}_*$ и осуществляем переход к нижнему пределу в \eqref{eqq-201}; с учетом оценки \eqref{eqq-401} получаем
$$
\liminf_{n\to\infty}\langle\mathbf{M}(\boldsymbol{v}_n),\boldsymbol{v}_*-\boldsymbol{v}_n\rangle\geq 0,
$$
откуда следует, что
$$
\limsup_{n\to\infty}\langle\mathbf{M}(\boldsymbol{v}_n),\boldsymbol{v}_n-\boldsymbol{v}_*\rangle\leq 0. 
$$

В силу псевдомонотонности оператора $\mathbf{M}$ приходим к неравенству:
\begin{equation}
\liminf_{n\to\infty}\langle\mathbf{M}(\boldsymbol{v}_n),\boldsymbol{v}_n-\boldsymbol{w}\rangle\geq \langle\mathbf{M}(\boldsymbol{v}_*),\boldsymbol{v}_*-\boldsymbol{w}\rangle\quad\forall\,\boldsymbol{w}\in \boldsymbol{V}.\label{eqq-601}
\end{equation}

Перепишем теперь неравенство \eqref{eqq-201} в более удобном для наших целей виде:
$$
\phi_g(\boldsymbol{v}_n)-\phi_g(\boldsymbol{w})\leq
-\langle\mathbf{M}(\boldsymbol{v}_n)+\mathbf{T}_{\boldsymbol{u}_n}(\boldsymbol{v}_n),\boldsymbol{v}_n-\boldsymbol{w}\rangle\quad \forall\,\boldsymbol{w}\in\boldsymbol{V},\;n\in\mathbb{N}. 
$$
Переходя в этом неравенстве к верхнему пределу (при $n\to\infty$), получим
$$
\limsup_{n\to \infty}\{\phi_g(\boldsymbol{v}_n)-\phi_g(\boldsymbol{w})\}\leq\limsup_{n\to \infty}\{-\langle\mathbf{M}(\boldsymbol{v}_n)+\mathbf{T}_{\boldsymbol{u}_n}(\boldsymbol{v}_n),\boldsymbol{v}_n-\boldsymbol{w}\rangle\}.
$$
Отсюда и из \eqref{eqq-401} и \eqref{eqq-601}, выводим неравенства
$$
\phi_g(\boldsymbol{v}_*)-\phi_g(\boldsymbol{w})\leq\liminf_{n\to \infty}\{\phi_g(\boldsymbol{v}_n)-\phi_g(\boldsymbol{w})\}\leq\limsup_{n\to \infty}\{\phi_g(\boldsymbol{v}_n)-\phi_g(\boldsymbol{w})\}
$$
$$
\leq\limsup_{n\to \infty}\{-\langle\mathbf{M}(\boldsymbol{v}_n)+\mathbf{T}_{\boldsymbol{u}_n}(\boldsymbol{v}_n),\boldsymbol{v}_n-\boldsymbol{w}\rangle\}=-\liminf_{n\to \infty}\{\langle\mathbf{M}(\boldsymbol{v}_n)+\mathbf{T}_{\boldsymbol{u}_n}(\boldsymbol{v}_n),\boldsymbol{v}_n-\boldsymbol{w}\rangle\}
$$
$$
\leq-\langle\mathbf{M}(\boldsymbol{v}_*)+\mathbf{T}_{\boldsymbol{u}_*}(\boldsymbol{v}_*),\boldsymbol{v}_*-\boldsymbol{w}\rangle\qquad\forall\,\boldsymbol{w}\in \boldsymbol{V}.
$$
и, следовательно,
$$
\langle\mathbf{M}(\boldsymbol{v}_*)+\mathbf{T}_{\boldsymbol{u}_*}(\boldsymbol{v}_*),\boldsymbol{w}-\boldsymbol{v}_*\rangle +\phi_g(\boldsymbol{w}) -\phi_g(\boldsymbol{v}_*)\geq 0\quad\forall\,\boldsymbol{w}\in \boldsymbol{V}.
$$
Тем самым доказано, что $(\boldsymbol{v}_*,\boldsymbol{u}_*)\in\boldsymbol{G}$. 

{\bf Шаг 4.} С учетом вышеизложенного мы имеем возможность применить лемму~\ref{lemm2}  для обоснования разрешимости задачи \eqref{bing-op-1}--\eqref{bing-op-7}, а равенство \eqref{est-gr} получается аналогично \eqref{eqq-20}. 

Таким образом, теорема~\ref{teor1} полностью доказана. 
\begin{flushright}
$\Box$
\end{flushright}

\newpage

\end{document}